\documentclass[11pt,twoside]{amsart}

\usepackage{enumerate}
\usepackage[hyperindex=true,plainpages=false,colorlinks=false,pdfpagelabels]{hyperref}
\theoremstyle{plain}
\newtheorem{The}{Theorem}[section]
\newtheorem*{The*}{Theorem}

\newtheorem{Lem}[The]{Lemma}
\newtheorem{Cor}[The]{Corollary}
\newtheorem*{Cor*}{Corollary}

\theoremstyle{definition}

\theoremstyle{remark}

\newtheorem*{Rem*}{Remark}

\numberwithin{equation}{section}

\DeclareMathOperator{\Span}{Span}

\DeclareMathOperator{\sff}{I\!I}

\renewcommand{\Im}{\operatorname{Im}}
\renewcommand{\Re}{\operatorname{Re}}

\newcommand{\zb}{{\bar z}}
\DeclareMathOperator{\dbar}{\bar\partial}

\newcommand{\R}{\mathbb{R}}
\newcommand{\C}{\mathbb{C}}

\DeclareMathOperator{\Will}{\mathcal{W}}

\begin{document}

\title[Constant mean curvature tori are Davey--Stewartson stationary]{Constant
  mean curvature tori as stationary solutions to the Davey--Stewartson
  equation}

\author{Christoph Bohle}

\address{Christoph Bohle\\
  Institut f\"ur Mathematik\\
  Technische Universit{\"a}t Berlin\\
  Stra{\ss}e des 17.\ Juni 136\\
  10623 Berlin\\
  Germany}

\email{bohle@math.tu-berlin.de}

%\keywords{constant mean curvature surface,Davey-Stewartson equation}


\date{\today}

\thanks{Author supported by DFG SPP 1154 ''Global Differential Geometry''.}

\begin{abstract} 
  A well known result of Da Rios and Levi--Civita says that a closed planar
  curve is elastic if and only if it is stationary under the localized
  induction (or smoke ring) equation, where stationary means that the
  evolution under the localized induction equation is by rigid motions.  We
  prove an analogous result for surfaces: an immersion of a torus into the
  conformal 3--sphere has constant mean curvature with respect to a space form
  subgeometry if and only if it is stationary under the Davey--Stewartson
  flow.
\end{abstract}

\maketitle

%%%%%%%%%%%%%%%%%%%%%%%%%%%%%%%%%%%%%%%%%%%%%%%%%%%%%%%%%%%%%%%%%%%%%%%%%%%%%%
%           Introduction                                                     %
%%%%%%%%%%%%%%%%%%%%%%%%%%%%%%%%%%%%%%%%%%%%%%%%%%%%%%%%%%%%%%%%%%%%%%%%%%%%%%

\section{Introduction}
\label{sec:intro}

The Davey--Stewartson (DS) hierarchy \cite{DS} is a 2+1--dimensional
generalization of the non--linear Schr\"odinger (NLS) hierarchy. The time
evolutions in both cases describe integrable deformations of Dirac operator
potentials in 1 or 2 dimensions, respectively.

It is well know since the work of Hasimoto \cite{Ha} that the NLS--hierarchy
has a geometric realization as an evolution of space curves describing the
motion of infinitely thin vortex filaments, the localized induction equation
\[ \dot \gamma = J \gamma''= \gamma'\times \gamma'' \] introduced by
Levi--Civita's student Da Rios \cite{DR} (see also~\cite{Ri}), where $J$
denotes the $90^\circ$--degree rotation in the normal bundle and $'$ the
derivative with respect to arc length: if a curve $\gamma$ in $\R^3$ moves
under the localized induction equation its complex curvature $\psi=\kappa
e^{i\int \tau}$ evolves under the non--linear Schr\"odinger equation. While
the localized induction equation describes the extrinsic evolution of curves
in space, the non--linear Schr\"odinger equation describes the corresponding
evolution of the differential invariants curvature $\kappa$ and
torsion~$\tau$. A curve is stationary under this smoke ring evolution if it
evolves by rigid motions and reparametrization only or, which is equivalent,
if its complex curvature is a traveling wave solution of the NLS--equation.

An analogous geometric version of the Davey--Stewartson flow as an evolution
of surfaces in 4--space was introduced by Konopelchenko \cite{Ko96,Ko00} who
realized that the Weierstrass representation, a correspondence between
immersions into Euclidean space and solutions to Dirac equations, allows to
deform conformal immersions of surfaces by deforming the corresponding Dirac
potentials.  While Konopelchenko's discussion is rather local, Taimanov
\cite{Ta97,Ta06} pointed out that the evolution equations can in fact be
globally defined for immersions of tori with trivial normal bundle and
preserves their conformal type.  A more geometric definition of the
Davey--Stewartson flow that avoids the use of Dirac operators and underlines
the analogy to the smoke ring evolution of space curves was given by Burstall,
Pedit, and Pinkall \cite{BPP02}: the normal part of the Davey--Stewartson
evolution of an immersion $f\colon T^2\rightarrow S^4$ into the conformal
4--sphere is
\[ \dot f^\perp = J \sff^\circ (X,X),\] where $J$ denotes $90^\circ$ rotation
in the normal bundle of $f$, $\sff^\circ$ is the trace free second fundamental
form with respect to some compatible space form subgeometry, and $X$ is some
holomorphic vector field. In case $f$ has topologically trivial normal bundle
there is an essentially unique way to complement $\dot f^\perp$ to a conformal
deformation $\dot f$ by adding a tangential deformation.  It should be noted
that, in contrast to the case of curves, the resulting conformal deformation
$\dot f$ of a torus $f\colon T^2\rightarrow S^4$ with trivial normal bundle is
not given explicitly, because the tangential deformation is obtained by
solving a $\dbar$--problem.  Unlike the smoke ring evolution of curves in
$\R^3$, the Davey--Stewartson flow depends on the additional choice of a
holomorphic vector field or, equivalently, the choice of conformal coordinates
on the universal covering of $T^2$; scaling these coordinates by some complex
factor changes the Davey--Stewartson flow. An immersion $f\colon
T^2\rightarrow S^4$ is called stationary under the Davey--Stewartson flow if,
for every choice of conformal coordinates, the evolution is by M\"obius
transformations and reparametrization only.

We prove the following analogue to Da Rios's and Levi--Civita's result that
the closed planar curves stationary under the localized induction equation are
plane elastica:
\begin{The*}
  A torus in the conformal 3--sphere is stationary under the Davey--Stewartson
  flow if and only if it is constrained Willmore and strongly isothermic.
\end{The*}

Here constrained Willmore means that $f\colon T^2\rightarrow S^3$ is a
critical point of the Willmore functional under infinitesimal conformal
variations and strongly isothermic means that $f$ is a critical point of the
projection to Teichm\"uller space, cf.~\cite{BPP1}.  A theorem of Richter
\cite{Ric97,BPP02} implies that a torus in $S^3$ is strongly isothermic and
constrained Willmore if and only if it has constant mean curvature with
respect to some space form subgeometry.

The main theorem of the paper thus provides a M\"obius geometric
characterization of constant mean curvature tori in space forms as stationary
solutions to the Davey--Stewartson equation.  In the appendix we prove a local
M\"obius geometric characterization of constant mean curvature surfaces in
space forms which is due to K.~Voss and uses Bryant's quartic
differential~\cite{Br84}.

\section{Invariants of Immersions into the Conformal
  $n$--Sphere}\label{sec:invariants}

We briefly review the Burstall--Pedit--Pinkall--approach \cite{BPP02} to
define invariants of immersions into the conformal $n$--sphere. For details
see Sections~3.1 to 3.3 of~\cite{BPP02}.

As a model for the \emph{conformal $n$--sphere} $S^n$ we use the projectivized
lightcone in Minkowski space $\R^{n+1,1}$ of dimension $n+2$, that is, in
$\R^{n+2}$ equipped with the metric \[<v,w> = -v_0 w_0+ v_1 w_1 + ... +
v_{n+1} w_{n+1}.\] Given an immersion $f\colon M \rightarrow S^n$ of a
manifold $M$ into the projective lightcone $S^n$, the pull back
$g=<d\psi,d\psi>$ of the Minkowski metric with respect to a homogeneous lift
$\psi$ of $f$ is a Riemannian metric on $M$ whose conformal class is
independent of $\psi$ because $\tilde g = \lambda^2 g$ for
$\tilde\psi=\lambda\psi$. The name conformal $n$--sphere is justified by the
fact that the standard metric of the round $n$--sphere in Euclidean space is
obtained from the embedding $x\in \{x \in \R^{n+1}\mid \sum x_i^2=1\} \mapsto
(1,x) \in S^n$ into the lightcone.  The identity component $O_o(n+1,1)$ of the
isometry group of $\R^{n+1,1}$ acts on the projective lightcone $S^n$ as the
group of orientation preserving M\"obius transformations.

Let $M$ now be a Riemann surface and $f\colon M \rightarrow S^n$ a
\emph{conformal immersion} which means that the conformal structure induced by
$f$ coincides with the Riemann surface structure of $M$.  The \emph{mean
  curvature sphere congruence} of $f$ is the sphere congruence spanned by
\[ \mathcal{V} = \Span\{ \psi, \psi_z,\psi_{\zb}, \psi_{z\zb}\}, \] where
$\psi$ is an arbitrary homogeneous lift of $f$ and $z$ are local conformal
coordinates on $M$.  One easily checks that $\mathcal{V}$ is a well defined
Minkowski space bundle of rank~4 and therefore indeed describes a congruence
of 2--sphere in $S^n$, namely the unique congruence of 2--spheres tangent to
$f$ for which the mean curvature at every point of contact coincides with that
of~$f$ (where the mean curvature is taken with respect to an arbitrary
compatible space form subgeometry).  The bundle $\mathcal{V}^\perp$ is called
the \emph{M\"obius normal bundle} and carries a positive definite metric and a
metric connection $D$, the \emph{normal connection}, defined by orthogonal
projection of the ordinary derivative on $\R^{n+1,1}$.

We assume now that $M$ is equipped with a fixed nowhere vanishing conformal
vector field $X$ or, equivalently, with a nowhere vanishing complex
holomorphic 1--form $dz$ dual to $X$.  Then there is a unique future pointing
homogeneous lift $\psi$ of $f$ normalized by the property that the metric
induced by $\psi$ coincides with $|dz|^2$.  We call this lift the
\emph{normalized lift} of $f$ with respect to~$X$.  While the lightcone model
of $S^n$ and the mean curvature sphere congruence date back to Darboux and
Thomsen, the normalized lift as an efficient means to define M\"obius
invariants of immersion was only introduced quite recently in \cite{BPP02}.

There is a unique section $\hat \psi$ of $\mathcal{V}$ that satisfying
$<\hat\psi,\hat\psi> =0$, $<\psi,\hat\psi> =-1$ and $<d\psi,\hat\psi> = 0$,
where $\psi$ is the normalized lift of $f$ with respect to~$z$. This yields a
normalized frame $(\psi, \psi_z,\psi_\zb ,\hat\psi)$ of $\mathcal{V}^\C$.
Since $\psi_{zz}$ is orthogonal to $\psi$, $\psi_z$ and $\psi_\zb$, there is a
complex function $c$ and a section $\kappa \in \Gamma((\mathcal{V}^\perp)^\C)$
such that
\begin{equation}
  \label{eq:hill_equiation}
  \psi_{zz} + \tfrac c2 \psi = \kappa.
\end{equation}
The quantities $\kappa$ and $c$ are called the \emph{conformal Hopf
  differential} and the \emph{Schwarzian derivative} of $f$.  Their dependence
on the choice of $dz$ is as follows: if $\tilde \kappa$ and $\tilde c$ are the
invariants with respect to $d\tilde z$, then
\begin{equation} \label{eq:change_of_kappa} \tilde \kappa \frac{d\tilde
    z^2}{|d\tilde z|} = \kappa \frac{dz^2}{|dz|} \qquad \textrm{ and } \qquad
  \tilde c d\tilde z^2 = (c-S_z(\tilde z)) dz^2,
\end{equation}
where $S_z(g)= \left(\frac{g''}{g'}\right)' - \frac12
\left(\frac{g''}{g'}\right)^2$ denotes the classical Schwarzian derivative of
a holomorphic function $g$ with respect to $z$.

Let $(\psi, \psi_z,\psi_\zb ,\hat\psi)$ be the frame of $\mathcal{V}^\C$ and let
$\xi\in \Gamma(\mathcal{V}^\perp)$ be a section of the M\"obius normal bundle. A
straightforward calculation shows that the frame equations are
\begin{equation}
  \label{eq:lc_frame}
  \begin{split}
  \psi_{zz} &= -\tfrac c2 \psi + \kappa \\
  \psi_{z\zb} &= -|\kappa|^2 \psi + \tfrac12 \hat \psi \\
  \hat\psi_z &= -2 |\kappa|^2 \psi_z -c \psi_\zb+ 2D_\zb \kappa \\
  \xi_z &= 2<D_\zb\kappa,\xi> \psi -2 < \kappa, \xi>\psi_\zb + D_z
  \xi.
  \end{split}
\end{equation}
The \emph{conformal Gauss, Codazzi} and \emph{Ricci equation} are the
integrability equations
\begin{gather}
  \label{eq:lc_gauss}
  \tfrac12 c_\zb = 2|\kappa|^2_z + <D_z\bar\kappa,\kappa> - <
  \bar\kappa,D_z\kappa> \\
  \label{eq:lc_Codazzi}
  \Im(D_\zb D_\zb\kappa+\tfrac {\bar c}2\kappa)=0 \\
  \label{eq:lc_Ricci}
  R^D_{z\zb}\xi = 2 <\bar\kappa,\xi> \kappa - 2< \kappa,\xi> \bar
  \kappa. 
\end{gather}
For immersions into the 4--sphere $S^4$ the normal bundle is 2--dimensional
and has a complex structure $J$ compatible with metric and orientation.  The
Ricci equation \eqref{eq:lc_Ricci} then reads
\begin{equation}
  R^D_{z\zb} = 2 <J \bar\kappa,\kappa > J \tag{\ref{eq:lc_Ricci}'}.
\end{equation}  
The invariants $\kappa$, $c$ and $D$ describe the immersions $f$ uniquely up
to M\"obius transformation. Conversely, given $\kappa$, $c$ and $D$ (with $D$
a metric connection on an abstract bundle $\mathcal{V}^\perp$ and $\kappa$ a
section of its complexification) that satisfy the conformal
Gauss--Codazzi--Ricci equations, there is a unique conformal immersions $f$
with M\"obius monodromy belonging to this data, cf.\ \cite{BPP02,Boh03}.

An immersion $f$ is called \emph{isothermic} if locally and away from umbilic
points it admits conformal curvature line coordinates.  One easily checks
(using (26) of \cite{BPP02}) that $z=x+iy$ are conformal curvature line
coordinates for an immersion $f$ if and only if the conformal Hopf
differential $\kappa$ of $f$ with respect to $z$ is a real section of the
M\"obius normal bundle~$\mathcal{V}^\perp$. Thus $f$ is isothermic if locally
and away from umbilics it admits conformal coordinates $z$ for which $\kappa$
is real.

If $M=T^2$ is a torus it makes sense to consider isothermic immersions $f$
admitting global conformal curvature line coordinates, that is, conformal
coordinates $z$ on the universal covering of $T^2=\C/\Gamma$ with respect to
which the conformal Hopf differential $\kappa$ is real.  Such immersions are
called \emph{strongly isothermic}.

It is proven in \cite{BPP02}, see equations (34) and (54) there, that a
conformal immersion of a torus $M=T^2$ is constrained Willmore, i.e., a
critical point of the Willmore energy $\Will=\int |\kappa|^2$ under
infinitesimal conformal variations, if and only if
\begin{equation}
  \label{eq:constrained_willmore_lightcone}
  D_\zb D_\zb \kappa + \tfrac12\bar c \kappa = \Re(\lambda \kappa )   
\end{equation}
for some $\lambda\in \C$, where $z$ are conformal coordinates on the universal
covering.

A theorem by J.~Richter \cite{Ric97} (proven in (36) to (40) of \cite{BPP02})
states that an immersion $f\colon T^2\rightarrow S^3$ of a torus into the
conformal 3--sphere is simultaneously strongly isothermic and constrained
Willmore if and only if it has constant mean curvature with respect to some
space form subgeometry.

It should be mentioned that if the underlying compact Riemann surfaces $M$ has
higher genus the definition of strongly isothermic immersions $f\colon M
\rightarrow S^n$ and the Euler--Lagrange equation for constrained Willmore
immersions $f\colon M \rightarrow S^n$ are more involved, because the
coordinates $z$ have to be replaced by a holomorphic quadratic differential,
see \cite{BPP1}.

\section{Davey--Stewartson Flow}\label{sec:flows}

We describe now the M\"obius geometric approach \cite{BPP02} to the
Davey--Stewartson flow on the space of immersions $f\colon T^2\rightarrow S^4$
with topologically trivial normal bundle. As a preparation we discuss the
effect of arbitrary infinitesimal conformal deformations of immersions
$f\colon T^2 \rightarrow S^4$ on the differential invariants of
Section~\ref{sec:invariants}. (For this we follow Section~4.1 of~\cite{BPP02}
and 14.3 of \cite{Boh03}.)

We fix conformal coordinates $z$ on the universal covering of the torus
$T^2$. An infinitesimal deformation of the homogeneous lift of a conformal
immersion $f\colon T^2 \rightarrow S^4$ is of the form
\begin{equation}
  \label{eq:deform}
  \dot \psi = a \psi + b \psi_z + \bar b \psi_\zb + \sigma
\end{equation}
with real function $a$, complex function $b$ and $\sigma \in
\Gamma(\mathcal{V}^\perp)$.  The deformation $\dot \psi$ hat no $\hat
\psi$--component, because we deform tangential to the lightcone. The
infinitesimal condition $< \dot \psi_z,\psi_z>=0$ for preserving the conformal
structure is equivalent to
\begin{equation}
  \label{eq:delbar}
  \bar b_z = 2<\sigma,\kappa>
\end{equation}
and the condition $\Re < \dot \psi_z,\psi_\zb>=0$ for preserving the
normalization is equivalent to
\[ a = - \Re b_z.\] The tangential part $b$ of a conformal deformation $\dot
f$ is thus essentially determined by the normal part $\sigma$ via the
$\dbar$--equation \eqref{eq:delbar} which can be solved if and only if
$\int<\sigma,\kappa>=0$.

In codimension~2, unlike in the codimension~1 case, the tangential and normal
deformations $b$ and $\sigma$ do not completely determine the deformation of
the Gauss--Codazzi--Ricci data because one can apply infinitesimal gauge
transformations of the normal bundle without changing $\psi$.  Such
infinitesimal normal bundle rotation is given by $\chi\colon T^2\rightarrow\R$
that describes the normal deformation of a section of the normal bundle $\xi
\in \Gamma(\mathcal{V}^{\perp})$ via
\begin{equation}
  \label{eq:nb_def}
  (\dot \xi)^\perp = \chi J \xi,
\end{equation}
where $J$ denotes the complex structure on the normal bundle.
    
The effect on the Gauss--Codazzi--Ricci--data $\kappa$, $c$ and $D$ of an
infinitesimal conformal deformation \eqref{eq:deform} together with an
infinitesimal normal bundle rotation \eqref{eq:nb_def} is
\begin{align}\label{eq:deformation_data}
  \dot \kappa & = D_zD_z \sigma + \tfrac c2 \sigma + (\tfrac32 b_z -\tfrac12
  \bar  b_\zb )\kappa + b D_z \kappa + \bar b D_\zb \kappa- \chi J\kappa  \\
  \dot c & = b_{zzz} + 2cb_z + bc_z + \bar bc_\zb + 16 < \sigma,\bar \kappa>
  <\kappa ,\kappa> +  \label{eq:deformation_data2}\\
  & \notag \qquad 6(< D_\zb D_z \sigma ,\kappa > - < D_z\sigma, D_\zb \kappa >
  ) + 2( < D_\zb \sigma ,D_z\kappa> - <\sigma ,D_\zb D_z \kappa >)
  \\
  \dot D_z \xi & = (\chi_z + 2 < J \sigma, D_\zb \kappa >+ 2 < J \kappa,\bar b
  \bar \kappa + D_\zb \sigma >) J \xi. \label{eq:deformation_data3}
\end{align}

We describe now the flows up to order two of the Davey--Stewartson hierarchy
on the space of immersions $f\colon T^2\rightarrow S^4$ with trivial normal
bundle.

The 0\textsuperscript{th}--order flow is given by $\sigma =0$, $b=0$ and
$\chi\colon T^2\rightarrow \R$ and hence consists of a pure normal bundle
rotation (the immersion does not actually move, but normal frames are
rotated).  By \eqref{eq:deformation_data} the resulting deformation of the
invariants is
\begin{align}
  \label{eq:rot_flow}
  \begin{split}
  \dot \kappa &= -\chi J\kappa \\
  \dot c &= 0 \\
  \dot D &=  d\chi J. 
  \end{split}
\end{align}

The 1\textsuperscript{st}--order flow is obtained by setting $\sigma =0$,
$b\in \C$ and $\chi =0$. This corresponds to a reparametrization of the
conformal immersion (without normal deformation) under which the invariants
evolve by
\begin{align}
  \label{eq:trans_flow}
  \begin{split}
  \dot \kappa &= b D_z \kappa + \bar b D_\zb \kappa \\
  \dot c &= b c_z + \bar b c_\zb \\
  \dot D_z &= 2 \bar b <J\kappa,\bar\kappa>J. 
  \end{split}
\end{align} 

The 2\textsuperscript{nd}--order flow, the \emph{Davey--Stewartson flow} with
respect to the coordinates $z$, is the conformal deformation whose normal part
is
\begin{equation}
  \label{eq:definition-DS}
  \sigma = 2\Re(J\kappa) = J\kappa+J\bar\kappa
\end{equation}
and whose tangential part $b$ satisfies the $\dbar$--equation
\eqref{eq:delbar}
\[ \bar b_z = 2< J\bar \kappa,\kappa>. \] Since we are on a torus such $b$
exists if and only if $\int_M < J\bar \kappa,\kappa>=0$ which, by the
Ricci--equation \eqref{eq:lc_Ricci}, is equivalent to the normal bundle degree
of the immersion being zero.  For tori $f\colon T^2\rightarrow S^4$ with
topologically trivial normal bundle the tangential deformation $b$ making the
Davey--Stewartson deformation conformal is thus uniquely defined up to adding
a constant. In other words, the 2\textsuperscript{nd}--order flow is then well
defined up to 1\textsuperscript{st}--order flow.  It should be noted that
changing the chosen coordinates $z$ changes the Davey--Stewartson flow.

\section{Davey--Stewartson Stationary Tori}

An immersed torus $f\colon T^2\rightarrow S^4$ is \emph{stationary under the
  Davey--Stewartson flow} if the Davey--Stewartson flows with respect to all
coordinates acts by M\"obius transformations and reparametrizations only.
Davey--Stewartson stationary tori generalize homogeneous tori, the tori
stationary under the 1\textsuperscript{st}--order ``reparametrization flow''
which are orbits of 2--parameter groups of M\"obius transformations and hence
M\"obius equivalent to tori of revolution whose profile curves are circles.

Due to the non--explicitness of the Davey--Stewartson equation caused by the
$\dbar$--equation describing the tangential part of the flow, the
investigation of general Davey--Stewartson stationary immersions $f\colon
T^2\rightarrow S^4$ into the 4--sphere is difficult (see
Lemma~\ref{lem:stationary_DS_condition} below for the equations characterizing
such tori).  In the following we focus on the special case of
Davey--Stewartson stationary tori that are contained in a
3--sphere~$S^3\subset S^4$ and hence have flat normal bundle such that the
$\dbar$--problem can be solved explicitly.

\begin{The} \label{Th:stationary_DS_in_S3} An conformal immersion $f\colon T^2
  \rightarrow S^3\subset S^4$ into the conformal 4--sphere that takes values
  in a 3--sphere is stationary under the Davey--Stewartson flow if and only if
  it is strongly isothermic and constrained Willmore.
\end{The}

Combined with Richter's theorem mentioned in Section~\ref{sec:invariants} this
implies:
\begin{Cor}\label{Cor:stationary_DS_in_S3}
  A conformal immersion of a torus into the conformal 3--sphere is stationary
  under the Davey--Stewartson flow if and only if it has constant mean
  curvature with respect to some space form subgeometry.
\end{Cor}

We derive now equations characterizing immersions $f\colon T^2 \rightarrow
S^4$ with topologically trivial normal bundle that are Davey--Stewartson
stationary.

\begin{Lem}
  \label{lem:stationary_DS_condition}
  A conformal immersion $f\colon T^2 \rightarrow S^4$ with trivial normal
  bundle is stationary under the Davey--Stewartson flow if and only if there
  are $\alpha$, $\beta\in \C$ and $\mu\colon T^2\rightarrow \C$ such that the
  invariants of $f$ with respect to a conformal chart $z$ on the covering of
  $T^2$ satisfy
\begin{gather}
  \label{eq:stationary_DS1}
  J (D_z D_z \kappa + \frac c2 \kappa) + \frac 32 b_z \kappa + b D_z
  \kappa = \alpha D_z \kappa + \bar \beta D_\zb \kappa - \mu J\kappa
  \\
  \label{eq:stationary_DS2}
  J (D_z D_z \bar\kappa + \frac c2 \bar
  \kappa) -
  \frac 12 \bar b_\zb \kappa + \bar b D_\zb \kappa = \beta D_z \kappa +
  \bar \alpha D_\zb \kappa - \bar \mu J\kappa \\
  \label{eq:stationary_DS3}
  \begin{split}
    b_{zzz} + 2c b_z& + b c_z + 16<J\kappa,\bar \kappa> <\kappa,\kappa> +
    \\ +& 8(<D_\zb D_z J\kappa,\kappa> - <D_z J\kappa, D_\zb \kappa>) =
    \alpha c_z + \bar \beta c_\zb
  \end{split} \\
  \label{eq:stationary_DS4}
  \begin{split}
    \bar b c_\zb + 6(<D_\zb& D_z J\bar \kappa,\kappa> - <D_z J\bar \kappa,
    D_\zb \kappa>)+\\ +& 2 (<D_\zb J\bar \kappa,D_z \kappa> - < J\bar
    \kappa, D_\zb D_z\kappa>)= \beta c_z + \bar \alpha c_\zb
  \end{split} \\
   \label{eq:stationary_DS5}
  0 = 2 \bar \beta<\bar \kappa , J \kappa> + (\mu)_z \\
  \label{eq:stationary_DS6}
  2<\bar \kappa,J\kappa> \bar b + 2< D_\zb \bar \kappa,\ \kappa> - 2<
  D_\zb \kappa,\bar \kappa> = 2 \bar \alpha <\bar \kappa , J \kappa> +
  (\bar \mu)_z,
\end{gather}
where $b\colon T^2\rightarrow \C$ is a solution to the $\dbar$--equation $\bar
b_z = 2< J\bar \kappa,\kappa>$.
\end{Lem}
\begin{proof}
  We fix conformal coordinates $z$ on the universal covering of the
  torus. Because being stationary under the Davey--Stewartson flow means that
  the second order flows with respect to all conformal coordinates on the
  universal covering depend linearly on the lower order flows, we express now
  the flows with respect to rotated coordinates $\tilde z = a z$ ($a\in \C_*$)
  in terms of the invariants with respect to $z$.

  Because real scalings of the coordinates correspond to real scalings of the
  flows it is sufficient to consider coordinates of the form $\tilde z= z
  e^{-i \theta}$ with $\theta \in \R$.  From \eqref{eq:change_of_kappa} we
  obtain that the invariants with respect to $\tilde z$ are $\tilde \kappa =
  \kappa e^{2i\theta}$ and $\tilde c = c e^{2i\theta}$. Furthermore, if $b$
  solves $\dbar$--problem for the tangential deformation with respect to $z$,
  then $\tilde b = b e^{i\theta}$ solves it with respect to $\tilde z$,
  because $\frac \partial{\partial \tilde z} = \frac \partial{\partial z} e^{i
    \theta}$.

  By \eqref{eq:deformation_data} to \eqref{eq:deformation_data3}, the
  evolution of the Gauss--Codazzi--Ricci data under the Davey--Stewartson flow
  with respect to $\tilde z$ with $\tilde b= b e^{i\theta}$ and $\chi=0$ is
\begin{align*}
  %\label{eq:DS_flow}
  \dot \kappa &= e^{2i\theta}( J (D_z D_z \kappa + \frac c2 \kappa) +
  \frac 32 b_z \kappa + b D_z \kappa ) + \\
  \notag & \qquad e^{-2i\theta}( J (D_z D_z \bar\kappa + \frac c2 \bar \kappa)
  - \frac 12 \bar b_\zb \kappa + \bar b D_\zb \kappa) \\
  \dot c &= e^{2i\theta}(b_{zzz} + 2c b_z + b c_z + 16<J\kappa,\bar \kappa>
  <\kappa,\kappa> + \\ &\qquad\qquad\qquad\qquad\qquad\qquad\qquad + 8(<D_\zb
  D_z J\kappa,\kappa> - <D_z J\kappa,
  D_\zb \kappa>)) + \\
  \notag & \quad\; e^{-2i\theta}( \bar b c_\zb + 6(<D_\zb D_z J\bar
  \kappa,\kappa> - <D_z J\bar \kappa, D_\zb \kappa>)+ \\&
  \qquad\qquad\qquad\qquad\qquad\qquad\qquad + 2 (<D_\zb J\bar
  \kappa,D_z \kappa> - < J\bar \kappa, D_\zb D_z\kappa>)) \\
  \dot D_z &= e^{-2i\theta}( 2<\bar \kappa,J\kappa> \bar b + 2 < D_\zb \bar
  \kappa,\ \kappa> - 2 < D_\zb \kappa,\bar \kappa> )J.
\end{align*}

Now $f$ is stationary if there is a complex function $\alpha(\theta)$
depending on $\theta$ and a real function $g(\theta,p)$ depending on $\theta$
and $p\in T^2$ such that, for every $\theta$, the above deformation equals
\begin{align*}
  \dot \kappa &= \alpha(\theta)e^{i \theta} D_z \kappa + \bar
  \alpha(\theta)e^{-i\theta} D_\zb
  \kappa- g J\kappa \\
  \dot c &= \alpha(\theta)e^{i\theta} c_z + \bar\alpha(\theta)e^{-i
    \theta}  c_\zb \\
  \dot D_z &= 2 \bar \alpha(\theta)e^{-i\theta} <J\kappa,\bar\kappa>J+
  g(\theta)_z J.
\end{align*} 

Fourier decomposition of $g$ and $\alpha$ immediately yields that the
condition for $f$ to be Davey--Stewartson stationary is
\eqref{eq:stationary_DS1} to \eqref{eq:stationary_DS6}.
\end{proof}

\begin{proof}[Proof (of Theorem \ref{Th:stationary_DS_in_S3})]
  Because $f$ takes values in a 3--sphere, its normal bundle is flat which by
  (\ref{eq:lc_Ricci}') is equivalent to $<J\bar \kappa,\kappa>=0$. Hence $b=0$
  solves the $\dbar$--equation $\bar b_z = 2< J\bar \kappa,\kappa>$ which
  drastically simplifies all equations in Lemma
  \ref{lem:stationary_DS_condition}.  In particular, since we are on a torus,
  equation \eqref{eq:stationary_DS5} implies that $\mu\colon T^2\rightarrow
  \C$ is constant.
  
  Because $f$ takes values in a 3--sphere there is a constant vector $n\in
  \R^{5,1}$ of length $1$ that is contained in every fiber of the M\"obius
  normal bundle $\mathcal{V}^\perp$ and satisfies
  \begin{equation}
    \label{eq:immersion_in_S3}
    <\kappa,n>=0.
  \end{equation}
  In particular $\kappa$ pointwise is a complex multiple of $Jn$. We use this
  in order to decompose equations \eqref{eq:stationary_DS1} and
  \eqref{eq:stationary_DS2} into $n$-- and $Jn$--parts
  \begin{align*}
    D_z D_z \kappa + \frac c2 \kappa = - \mu \kappa
    \tag{\ref{eq:stationary_DS1}'} \\
    D_z D_z \bar\kappa + \frac c2 \bar \kappa = - \bar \mu \kappa
    \tag{\ref{eq:stationary_DS2}'}
  \end{align*}
  and
  \begin{align*}
    \alpha D_z \kappa + \bar \beta D_\zb \kappa =0
    \tag{\ref{eq:stationary_DS1}''} \\
    \beta D_z \kappa + \bar \alpha D_\zb \kappa
    =0\tag{\ref{eq:stationary_DS2}''}.
  \end{align*}
  Moreover, the remaining scalar product terms in \eqref{eq:stationary_DS3}
  and \eqref{eq:stationary_DS4} vanish such that
  \begin{align*}
     \alpha c_z + \bar \beta c_\zb=0 \tag{\ref{eq:stationary_DS3}'}\\
     \beta c_z + \bar \alpha c_\zb=0.\tag{\ref{eq:stationary_DS4}'}
  \end{align*}
  The equations (\ref{eq:stationary_DS1}''), (\ref{eq:stationary_DS2}''),
  (\ref{eq:stationary_DS3}') and (\ref{eq:stationary_DS4}') can only be solved
  by $\alpha= \beta=0$ (unless $f$ is ``equivariant'' with respect to a
  1--parameter group of M\"obius transformations).
  
  If $\mu \neq 0$ the Codazzi equation \eqref{eq:lc_Codazzi} implies that the
  imaginary part of (\ref{eq:stationary_DS2}') is $\Im(\bar \mu
  \kappa)=0$. After rotation of the coordinates $z$ we can thus assume that
  $\kappa$ and $\mu$ are real. This shows that $f$ is strongly isothermic. On
  the other hand, taking the real part of (\ref{eq:stationary_DS2}') implies
  \[ D_\zb D_\zb \kappa + \frac {\bar c}2 \kappa = \Re(- \mu \kappa) = - \mu
  \kappa \] which shows that $f$ is constrained Willmore, see
  \eqref{eq:constrained_willmore_lightcone}.
  
  If $\mu=0$, then $f$ is Willmore. Because $f$ takes values in the 3--sphere,
  \eqref{eq:stationary_DS6} takes the form
  \[ < D_\zb \bar \kappa,\ \kappa> - < D_\zb \kappa,\bar \kappa> = 0
  \tag{\ref{eq:stationary_DS6}' }\] which, after writing $\kappa$ as
  $\kappa=a Jn$ with $a\colon T^2\rightarrow \C$, becomes $\bar a_\zb a -
  a_\zb \bar a=0$. Away from umbilic points the ratio $\bar a/a$ is constant
  such that the argument of $a$ is a locally constant function defined on the
  complement of the set of umbilic points. We can thus assume, possibly after
  rotation of $z$, that $\kappa$ is real on a component of the complement of
  the set of umbilic points. Now Willmore surfaces are analytic because their
  Euler--Lagrange equation is elliptic such that imaginary part of $\kappa$
  has to vanish globally and $f$ is strongly isothermic.

  Reversing the argumentation one immediately shows that a strongly
  isothermic, constrained Willmore immersion $f\colon T^2\rightarrow S^3$ is
  Davey--Stewartson stationary: choosing $z$ such that $\kappa$ is real,
  \eqref{eq:constrained_willmore_lightcone} together with the Codazzi equation
  implies $D_\zb D_\zb \kappa + \tfrac12 \bar c \kappa = \lambda \kappa$ with
  $\lambda\in \R$.  Thus, \eqref{eq:stationary_DS1} to
  \eqref{eq:stationary_DS6} are satisfied with $\alpha=\beta=b=0$ and
  $\mu=-\lambda$.
\end{proof}

\section*{Appendix}

We prove a (local) characterization of Willmore surfaces and constant mean
curvature surfaces in space forms in terms of Bryant's quartic
differential~\cite{Br84}.  This characterization was apparently first obtained
by K.~Voss~\cite{Vo} but seems to be nowhere published.  Our proof uses the
methods of Burstall, Pedit, and Pinkall~\cite{BPP02} which are described in
Section~\ref{sec:invariants}.

For surfaces in the conformal 3--sphere, the M\"obius normal bundle
$\mathcal{V}^\perp$ has a canonical trivialization by the unique section $Y\in
\Gamma(\mathcal{V}^\perp)$ of unit length compatible with the
orientation. Using this trivialization, the conformal Hopf differential with
respect to a chart~$z$ becomes a complex function $\kappa$ and the Gauss and
Codazzi equations \eqref{eq:lc_gauss} and \eqref{eq:lc_Codazzi} become
\[ \tfrac12 c_{\bar z} = (|\kappa|^2)_{z} + 2 \bar\kappa_z\kappa \qquad
\textrm{ and } \qquad \Im(\kappa_{\bar z\bar z} + \tfrac{\bar c}2 \kappa )=0.
\]
Equation \eqref{eq:lc_frame} shows that in our setting Bryant's quartic
differential $\mathcal{Q} = <Y_{zz},Y_{zz}> dz^4$ takes the form
\[ \mathcal{Q} = 4(\kappa \kappa_{\bar z z} + |\kappa|^2\kappa^2-\kappa_{\bar
  z}\kappa_z) dz^4. \] Holomorphicity of $\mathcal{Q}$ is equivalent to
\[ \kappa \kappa_{\bar z \bar z z} + (|\kappa|^2)_{\bar z}\kappa^2 +
2|\kappa|^2 \kappa \kappa_{\bar z} -\kappa_{\bar z\bar z}\kappa_z = 0.\] The
Gauss--equation implies that, away from umbilic points,
\[ \kappa^2\left(\frac{\kappa_{\bar z\bar z} + \tfrac{\bar c}2
    \kappa}\kappa\right)_z= \kappa \kappa_{\bar z \bar z z} +
(|\kappa|^2)_{\bar z}\kappa^2 + 2|\kappa|^2 \kappa \kappa_{\bar z}
-\kappa_{\bar z\bar z}\kappa_z.\] Thus, away from umbilic points,
holomorphicity of $\mathcal{Q}$ is equivalent to
\[ \kappa_{\bar z\bar z} + \tfrac{\bar c}2 \kappa = \lambda \kappa
\tag{$*$} \] for some anti--holomorphic function $\lambda$. The function
$\lambda$ vanishes identically if and only if the immersion is Willmore. In
case $\mathcal{Q}$ is holomorphic and $\lambda$ does not vanish identically,
away from the isolated zeroes of $\lambda$ we can introduce local coordinates
$\tilde z$ for which $\frac{d\tilde z}{dz}=\sqrt{\bar \lambda}$. By
\eqref{eq:change_of_kappa} we then have $\tilde \kappa |\lambda|^{3/2}= \kappa
\lambda$ such that $\tilde \kappa$ is real, because $\Im(\lambda \kappa)=0$ by
($*$) and the Codazzi equation. Without loss of generality we can thus assume
that $\kappa$ is real for the original coordinates $z$. The Codazzi equation
then shows that ($*$) holds with constant $\lambda$. But this implies that the
immersion has constant mean curvature with respect to some space form
subgeometry, see the paragraph in \cite{BPP02} that starts with equation (39).

This proves:
\begin{The}[K.Voss]
  The quartic differential $\mathcal{Q}$ of an immersion into the conformal
  3--sphere is holomorphic if and only if, away from umbilics and isolated
  points, the immersion is Willmore or has constant mean curvature with
  respect to some space form subgeometry.
\end{The}

%%%%%%%%%%%%%%%%%%%%%%%%%%%%%%%%%%%%%%%%%%%%%%%%%%%%%%%%%%%%%%%%%%%%%%%%%%%%%%
%           Bibliography                                                     %
%%%%%%%%%%%%%%%%%%%%%%%%%%%%%%%%%%%%%%%%%%%%%%%%%%%%%%%%%%%%%%%%%%%%%%%%%%%%%%

\end{document}